\documentclass[a4paper,twoside]{article}
\usepackage{RR}
\usepackage{hyperref}
\usepackage{appendix}
\newcommand{\lo}[1]{\|#1\|_0}
\def\A{\mathbf{A}}
\def\tA {\tilde{\A}}
\def\wA {\widehat{\A}}
\def\a {a}
\def\ta {\tilde{\a}}
\def\wa {\widehat{\a}}

\def\R{\mathbb{R}}
\def\C{\mathbb{C}}
\def\F{\mathbf{F}}

\def\lzero{\ell^0}

\def\lp{\ell^p}

\newenvironment{proof}[1][Proof]{\begin{trivlist}
\item[\hskip \labelsep {\bfseries #1}]}{\end{trivlist}}

\newcommand{\nrofsrc}{N} 
\newcommand{\nrofchn}{M}
\newcommand{\filterlength}{L}

\newtheorem{lemma}{Lemma}
\newtheorem{coro}{Corollary}
\newtheorem{theorem}{Theorem}
\newtheorem{definition}{Definition}

\usepackage[latin1]{inputenc}
\usepackage{graphicx,calc}

\usepackage{amssymb,amsmath}
\usepackage{stmaryrd}
\input cyracc.def
\font\tencyr=wncyr10
\def\cyr{\tencyr\cyracc}
\def\dirac{\mbox{\cyr sh}}

\newcommand{\qed}{\hfill \mbox{\raggedright \rule{.07in}{.1in}}}

\RRtitle{Caract\`ere bien pos\'e du probl\`eme de permutation pour l'estimation des filtres parcimonieux par minimisation $\lp$}
\RRetitle{Well-posedness of the permutation problem in sparse filter estimation with $\lp$ minimization\ \ }
\RRmotcle{Filtres parcimonieux, s\'eparation aveugle de sources, m\'elange convolutif, ambigu\"it\'e de permutation, th\'eor\`eme de mariage de Hall, matrice bi-stochastique}
\RRkeyword{sparse filter, convolutive blind source separation,  permutation ambiguity, $\lp$ minimization, Hall's Marriage Theorem, bi-stochastic matrix
}
\RRdate{November 2011}
\RRabstract{
Convolutive source separation is often done in two stages: 1) estimation of the mixing filters and 2) estimation of the sources. Traditional approaches suffer from the ambiguities of arbitrary permutations and scaling in each frequency bin of the estimated filters and/or the sources, and they are usually corrected by taking into account some special properties of the filters/sources. This paper focusses on the filter permutation problem in the absence of scaling, investigating the possible use of the temporal sparsity of the filters as a property enabling permutation correction. Theoretical and experimental results highlight the potential as well as the limits of sparsity as an hypothesis to obtain a well-posed permutation problem.
}
\RRresume{La s\'eparation de source des m\'elanges convolutifs se fait souvent en deux \'etapes : 1) estimation des filtres de m\'elange et 2) estimation des sources. Les approches classiques souffrent d'ambigu\"it\'es de permutation et de facteur d'\'echelle arbitraire pour chaque fr\'equence des filtres et/ou des sources estimés. Ces ambigu\"it\'es sont habituellement corrig\'ees en prenant en compte des propri\'et\'es particuli\`eres des filtres/sources. Cet article se concentre sur le probl\`eme de permutation des filtres en l'absence de facteur d'\'echelle, en explorant l'utilisation potentielle de la parcimonie temporelle des filtres pour r\'esoudre le probl\`eme de permutation. Les r\'esultats th\'eoriques et exp\'erimentaux soulignent tant le potentiel que les limites de l'hypoth\`ese de parcimonie pour obtenir un probl\`eme bien pos\'e.}
\RRprojet{METISS}
\RCRennes
\RRauthor{ Alexis~Benichoux,
Prasad~Sudhakar,
Frédéric~Bimbot,~and~R\'emi~Gribonval}

\begin{document}
\RRNo{7782}
\makeRR

\section{Introduction}
\label{sec:introduction}
%
%
%
%
Blind source separation and blind source localization are ubiquitous problems in signal processing, with applications ranging from wireless telecommunications to underwater acoustics and sound enhancement.

These problems can be considered as reasonably well understood and solved in simple linear instantaneous settings, where tools such as Independent Component Analysis, as well as techniques exploiting source sparsity, are now mature. However, the convolutive source localization / separation problem remains much more challenging. In particular, without further assumption than statistical independence between sources, the problem is known to be ill-posed because of the so-called frequency permutation (and scaling) problem: at best, one can hope to estimate for each frequency (up to a source and frequency dependent scaling factor) the collection of frequency components of all  sources (and of the associated mixing filters); but one cannot match the estimated frequency components from different subbands to globally identify the sources (and mixing filters). 

Several practical approaches have been proposed to solve the permutation and scaling problems in practice, by exploiting various properties of either the mixing filters or the sources to match different frequency subbands. While some of these methods may succeed in practice for certain types of sources / filters, there is no known theory guaranteeing the well-posedness of the permutation and scaling problem under appropriate assumptions. 

This paper contributes to fill this gap, by providing well-posedness guarantees for the permutation problem  under {\em sparsity} assumptions on the mixing filters. Sparse filters,  associated to impulse responses corresponding to a limited set of echoes, are typically encountered in a number of underwater communication channels~\cite{5352256} or wireless telecommunications scenarios~\cite{4518398, 5454399} which are relevant for blind source localization and separation, and the theoretical results achieved in this paper indicate that this property can potentially be exploited for blind estimation in this context.

\subsection{Problem formulation and notations}

Let
 $x_{i}[t]$, $1 \leq i \leq \nrofchn$ be $\nrofchn$ mixtures of $\nrofsrc$ source signals $s_j[t]$, resulting from the convolution with filters $\a_{ij}[t]$ of length $\filterlength$ such that:
\begin{equation}
x_i[t] = \sum_{j=1}^{\nrofsrc}(\a_{ij}\star s_j)[t],\quad 1\leq i\leq \nrofchn,
\label{eqn:timedomainmixing} 
\end{equation}
\noindent
where  $\star$ denotes convolution. The filter $\a_{ij}[t]$ typically models the impulse response between the $j^{th}$ source and the $i^{th}$ sensor.
By abuse of notation, $\F a_{ij} = \{a_{ij}[\omega]\}_{0\leq \omega <\filterlength}$ denotes the discrete Fourier transform of the filter seen as a vector $\a_{ij} = \{\a_{ij}[t]\}_{0\leq t < \filterlength} \in \C^{\filterlength}$. Also, the mixing equation~\eqref{eqn:timedomainmixing} can be rewritten as $\mathbf X=\mathbf A \star \mathbf S$, with $\A$ the matrix of filters 
\[
\A:=\left(\{\a_{ij}[t]\}_{0\leq t < \filterlength }\right)_{1\leq i\leq \nrofchn,\  1\leq j\leq \nrofsrc},
\]
$\mathbf X$ the observation matrix and $\mathbf S$ the source matrix.

In this context, blind filter estimation refers to the problem of obtaining estimates of the filters $\A$ from the mixtures $\mathbf X$, without any explicit knowledge about the sources $\mathbf S$. Mixing filters estimation is relevant for several purposes such as deconvolution, source localization, etc.~\cite{benestybook}. It also has a relationship with the problem of Multiple-Input-Multiple-Output (MIMO) system identification in communications engineering~\cite{934137}.

\subsection{Frequency domain filter estimation}
Estimating the mixing parameters is made easier when all filters are instantaneous, that is to say of length $\filterlength=1$, as the convolution product in~\eqref{eqn:timedomainmixing} is replaced by the usual product. However, things get complicated in the general setting of convolutive mixtures.  

A common approach for filter estimation then relies on the transformation of the mixing model in Eq.~\eqref{eqn:timedomainmixing} into the time-frequency domain, converting a single convolutive filter estimation problem into several complex instantaneous filter estimation problems. Using standard techniques for instantaneous mixing parameter estimation \cite{2010_AcademicPress_BookBSS}, complex mixing filter coefficients $$\tA[\omega]=\{\tilde {a}_{ij}[\omega]\}_{1\leq i\leq \nrofchn,\ 1\leq j\leq \nrofsrc}$$ are estimated for each frequency bin $0\leq \omega < \filterlength$. 


\subsection{Permutation and scaling ambiguities}  
Without further assumption on either the filters $\a_{ij}[t]$ or the sources $s_j[t]$, one can at best hope to find an estimation $\tA=(\ta_{ij})$ where for every frequency $\omega$ we have
\begin{equation}
\label{eq:DefPermScal}
\ta_{ij}[\omega]=  \lambda_j[\omega]a_{i\sigma_\omega(j)}[\omega],
\end{equation} with $\lambda_j[\omega]$ a scaling ambiguity and $\sigma_\omega \in {\mathfrak{S}}_N$ a permutation ambiguity, where $\mathfrak S_N$ is the set of permutations of the integers between one and $N$. 
Several methods  \cite{97} attempt to solve for these ambiguities by exploiting properties of either the sources $\mathbf S$ or the filters $\mathbf A$~\cite{serviere2004, 1327173, wang:a}.

\subsection{Exploiting sparsity to solve the permutation ambiguity}
The focus of this article is the use of the {\em sparsity} of $\mathbf A$ in the time domain to find $\sigma_0 \ldots \sigma_{\filterlength-1} \in \mathfrak S_N $, assuming the scaling is solved, i.e., $\lambda_{j}[\omega]=1$. 

Assuming that $\A$ is sparse means that each filter $\a_{ij}$ has  few nonzero coefficients, as measured by the $\lzero$ pseudo-norm
\[
\lo{\a_{ij}} := \sharp \{0 \leq t < \filterlength,\ \a_{ij}[t] \neq 0\} = \sum_{t} |a_{ij}[t]|^{0}.
\]

The approach considered in this article is to to seek permutations $\widehat\sigma_{0},\ldots \widehat\sigma_{\filterlength-1}$ yielding the sparsest estimated time-domain matrix of filters $\wA = (\wa_{ij})$ 
where $\wa_{ij}[\omega] := \ta_{i\widehat{\sigma}_{\omega}(j)}[\omega]$. Besides the $\lzero$ pseudo-norm $\lo{\wA} := \sum_{ij} \lo{\wa_{ij}}$, the following $\lp$ quasi-norms will be used
to quantify the sparsity of $\wA$:
\[
\|\wA\|_{p}^{p}:= \sum_{ij} \|\wa_{ij}\|_{p}^{p} = \sum_{ijt} |\wa_{ij}[t]|^{p},\quad 0 < p \leq 1.
\]

\subsection{Main results}
Our main result (Theorem~\ref{th:main}) is a theoretical guarantee that when the filter length $L$ is prime, $k$-sparse filters (i.e., such that $\|\a_{ij}\|_{0} \leq k$) uniquely minimize the $\lzero$ norm of $\A$ (up to a global permutation) if $\frac k L \leq \alpha(N)$, where $N$ is the number of sources. To reach this bound we exploit uncertainty principles as well as the bistochastic structure of the problem through an apparently new quantitative result on bistochastic matrices (Lemma~\ref{le:BiStochastic}).


\subsection{Structure of the paper}

The main theorems are stated in Section~\ref{sec:theory}, and the main ingredients of their proofs are described in Section~\ref{sec:ingredients}. In Section~\ref{sec:discussion} we discuss the strength of the assumptions used in the theorems, and how much these could be relaxed.
In Sec.~\ref{sec:algorithm}, a naive combinatorial $\lp$ minimization algorithm is proposed to resolve filter permutations and used for Monte-Carlo simulations. We conclude with a discussion of the potential, as well as the limits, of sparsity as a hypothesis to solve permutation problems, in connection with the theoretical and empirical results. All proofs are gathered in the appendix.

\section{Theoretical guarantees}
\label{sec:theory}
Given an $\nrofchn \times \nrofsrc$ filter matrix $\A$, made of filters of length $\filterlength$, and an $L$-tuple $(\sigma_{0},\ldots \sigma_{\filterlength-1}) \in \mathfrak{S}_N$ of permutations, we let $\tA$ be the matrix obtained from $\A$ by applying the permutations in the frequency domain, as in~\eqref{eq:DefPermScal}, without scaling ($\lambda_{j}[\omega]=1$). 

The effect of the permutations is said to coincide with that of a global permutation $\pi \in \mathfrak{S}_{N}$ of the columns of $\A$ if $\ta_{ij} = \a_{i\pi(j)}$, $\forall i,j$, or equivalently in the frequency domain:
\[
\ta_{ij}[\omega] := \a_{i\sigma_\omega(j)}[\omega] = \a_{i\pi(j)}[\omega],\ 0 \leq \omega < \filterlength,\ \forall i,j.
\]
This is denoted $\A \equiv \tA$. First, we show that for filters with disjoint time-domain supports, permutations cannot decrease the $\lp$ norm, $0 \leq p \leq 1$:

\begin{theorem}\label{th:disj}
Let $\Gamma_{ij} \subset\{0,\ldots,\filterlength-1\}$ be the time-domain support of $\a_{ij}$. Suppose that for all $i$ and $j_1\neq j_2$ we have 
\begin{equation}
\label{eq:disjointsupports}
\Gamma_{i,j_1} \cap \Gamma_{i,j_2} = \emptyset.
\end{equation}
Then, for $0\leq p\leq 1$, we have $\|\tA\|_p\geq\|\A\|_p$.
\end{theorem}

Note that filters with disjoint supports need not be very sparse: $\nrofchn$ filters of length $\filterlength$ can have disjoint supports provided that $\max_{j}\lo{\a_{ij}} \leq \filterlength/\nrofchn$. 
Yet, disjointness of filter supports is a strong assumption, and Theorem~1 only indicates that frequency permutations cannot decrease the $\lp$ norm. Thus, the minimum value of the $\lp$ norm might not be uniquely achieved (up to a global permutation). In our main result, we consider $k$-sparse filters of prime length, and $p=0$:
\begin{theorem}\label{th:main}
Let $\A$ be an $\nrofchn \times \nrofsrc$ matrix of filters of prime length $L$. Assume that 
\begin{equation}
\label{eq:ksparsefilters}
\max_{ij}\|{\a_{ij}}\|_0\leq k,
\end{equation} 
where
\begin{equation}
\label{eq:sparsityNsrc}
\frac kL \leq \alpha(N) := \left\{
\begin{array}{ll}
\frac{2}{N(N+2)} & \mbox{if $N$ is even},\\
\frac{2}{(N+1)^{2}} & \mbox{if $N$ is odd}.
\end{array}
\right.
\end{equation}
Then, up to a global permutation, $\A$ uniquely minimises the $\lzero$ pseudo-norm among all possible frequency permutations.   
\end{theorem}

\section{Main elements of the proof of Theorem~\ref{th:main}}
\label{sec:ingredients}
The proof of Theorem~\ref{th:main} relies on a measure of the ``amount'' of incurred permutation, on uncertainty principles, and on combinatorial arguments related to bi-stochastic matrices, involving Hall's Marriage Theorem.

\subsection{Measures of the amount of incurred permutations}
To measure the ``amount'' of incurred permutation, one can count the number of frequency bands where a non-trivial permutation is incurred, with respect to the best matching reference global permutation $\pi$, i.e., $\min_{\pi} \sharp \{\omega, \sigma_{\omega} \neq \pi\}$. However, this generally yields the maximum count $L-1$.

An alternative is to count the ``size'' of the incurred permutations,  given a reference global permutation $\pi$,  as the maximum number of frequencies where each estimated filter actually differs from the (globally permuted) original filters, yielding:
\begin{align}
\Delta(\tA,\A|\pi)
&:=
\max_{i,j}\lo{\F (\ta_{ij}-a_{i\pi(j)})}\\
\Delta(\tA,\A) & := \min_{\pi\in\mathfrak S_N} \Delta(\tA,\A|\pi).
\end{align}
Note that 
$\Delta(\tA,\A)=0$  iff $\tA \equiv \A$.

\subsection{Exploitation of an uncertainty principle}
With this notation, we have the following Lemma:
\begin{lemma}\label{le:uncertainty}
Assume that $\tA \not\equiv \A$, that $L$ is a prime integer, and that~\eqref{eq:ksparsefilters} holds with
\begin{equation}
\label{eq:sparsityassumption0prime}
2k + \Delta \leq L.
\end{equation}
Then $\lo{\tA} > \lo{\A}$ and
\(
\lo{\ta_{ij}} \geq \lo{\a_{ij}}, \forall i,j.
\)
The latter inequality is strict when $\ta_{ij} \neq \a_{ij}$. For a general $L$ (not necessarily prime), the same conclusions hold when the assumption~\eqref{eq:sparsityassumption0prime} is replaced with
\begin{equation}
\label{eq:sparsityassumption0}
2k \cdot \Delta < L.
\end{equation}
\end{lemma}

The skilled reader will rightly sense the role of uncertainty principles~\cite[Theorem 1]{Donoho:1989aa,elad2002generalized,Tao:2005aa} in the above lemma. 

\subsection{Combinatorial arguments}
Using Lemma~\ref{le:uncertainty} with prime $L$, a simple combinatorial argument can be used to obtain a weakened version of Theorem~\ref{th:main}, with the more conservative constant $\alpha'(N) := 1/2N!$:
by the pigeonhole principle, for any $L$-tuple of frequency permutations among $\nrofsrc$ sources, at least $L/\nrofsrc!$ permutations are identical; as a result, $\Delta(\tA,\A)$ is universally bounded from above by $L-L/\nrofsrc!$; hence if $k \leq L/2N!$ we obtain $2k +\Delta \leq L$ and we can conclude thanks to Lemma~\ref{le:uncertainty}.

The proof of Theorem~\ref{th:main} with the constant $\alpha(N)$ exploits a stronger universal upper bound $\Delta(\tA,\A) \leq L(1-2\alpha(N))$, obtained through an apparently new quantitative application of Hall's Marriage Theorem~\cite{Hall:1935aa} to bi-stochastic matrices.
\begin{definition}[Bi-stochastic matrix]
An $N \times N$ matrix $\mathbf{B}$ is called {\em bi-stochastic} if all its entries are non-negative, and the sum of the entries over each row as well as the sum of the entries over each column is one.
\end{definition}

\begin{lemma}
\label{le:BiStochastic}
Let $\mathbf{B}$ be an $N \times N$ bi-stochastic matrix: there exists a permutation matrix $\mathbf{P}$ such that all the entries of $\mathbf{B}$ on the support of $\mathbf{P}$ exceed the threshold
\begin{equation}
\label{eq:BiStochasticThreshold}
2\alpha(N) = \left\{
\begin{array}{ll}
\frac{4}{N(N+2)} & \mbox{if $N$ is even},\\
\frac{4}{(N+1)^{2}} & \mbox{if $N$ is odd}.
\end{array}
\right.
\end{equation}
\end{lemma}
\begin{coro}
\label{cor:perms}
Let $\sigma_{0}, \ldots, \sigma_{L-1} \in \mathfrak{S}_N$ be $L$ permutations. There exists a global permutation $\pi$ such that
\[
C_{j\pi(j)} = \sharp \{\ell: \sigma_{\ell}(j) = \pi(j)\} \geq 2L \alpha(N),\quad \forall 1 \leq j \leq N.
\]
\end{coro}


\section{Discussion}
\label{sec:discussion}
The reader may have noticed that Theorem~\ref{th:main}, while dropping the {\em disjoint support} assumption from Theorem~\ref{th:disj}, introduces new restrictions: the assumption that $L$ is prime, and the restriction to $p=0$ compared to $0\leq p \leq 1$ in Theorem~\ref{th:disj}.
How important are these restrictions ? Could they be relaxed while exploiting sparsity together with the {\em disjoint support} assumption ? This is discussed in this  section.

\subsection{Extending Theorem~\ref{th:main} to non-prime filter length $L$?}
As indicated by Lemma~\ref{le:sharpuncertainty} below, for even $L\geq 4$, there exists  sparse matrices of filters that are the sparsest {\em but not unique (even up to a global permutation)} solution of the considered problem: certain frequency permutations provide an {\em equally sparse but not equivalent solution}. 
\begin{lemma}\label{le:sharpuncertainty}
For any integer $k$ such that $2k$ divides $L$, there exists a matrix of $k$-sparse filters $\A$ and a set of $L/2k$ frequency permutations resulting in $\tA \not\equiv \A$ such that
for all $0 \leq p \leq \infty$: $\|\tA\|_{p}=\|\A\|_{p}$, and 
\begin{equation}
\|\ta_{ij}\|_{p}=\|\a_{ij}\|_{p},\quad \forall i,j.
\end{equation}
We have $2k \cdot \Delta(\tA,\A) =L$.
\end{lemma}

The fact that the filter matrices $\A$ and $\tA$ satisfy $2k \cdot \Delta(\tA,\A) =L$ shows the sharpness of Lemma~\ref{le:uncertainty} for the case when $L$ is even: the strict inequality in~\eqref{eq:sparsityassumption0} cannot be improved. 

Specializing Lemma~\ref{le:sharpuncertainty} to $k=1$ for even $L\geq 4$ yields ideally $1$-sparse filters $\a_{ij}$ and a set of $L/2$ frequency permutations such that: $\ta_{ij}$ are $1$-sparse; $\tA$ is not equivalent to $\A$ and cannot be discriminated from it by any $\ell^{p}$ norm. 

\subsection{Stronger guarantees with disjoint supports {\em and} sparsity ?}
Could one get improved results by combining the disjoint support assumptions from  Theorem~\ref{th:disj} and the sparsity assumption from Theorem~\ref{th:main} ? 
For even $L\geq 4$, Lemma~\ref{le:sharpuncertaintydisj} below indicates the existence of  sparse matrices of filters {\em with disjoint supports} that are the sparsest but not unique (even up to a global permutation) solution of the considered problem: certain frequency permutations of ``size'' $\Delta = L/2k$ provide an equally good but not equivalent solution. 
\begin{lemma}\label{le:sharpuncertaintydisj}
For any integers $k'<k\leq L/2$ such that $2k'$ divides $L$, there exist a matrix of $k$-sparse filters $\A$ {\em with disjoint supports}~\eqref{eq:disjointsupports}, and a set of $L/2k'$ frequency permutations resulting in $\tA \not\equiv \A$, such that for all $0 \leq p \leq \infty$: $\|\tA\|_{p}=\|\A\|_{p}$ and 
\begin{equation}
\|\ta_{ij}\|_{p}=\|\a_{ij}\|_{p},\quad \forall i,j.
\end{equation}
We have $2k' \cdot \Delta(\tA,\A)= L$.
\end{lemma}

Specializing Lemma~\ref{le:sharpuncertaintydisj} to $k'=1$ and $k=2$ for even $L\geq 4$ yields $2$-sparse filters $\a_{ij}$ and a set of $L/2$ frequency permutations such that: $\ta_{ij}$ are $2$-sparse; $\tA$ is not equivalent to $\A$ and cannot be discriminated from it by any $\ell^{p}$ norm.

This shows that even by adding the disjoint support assumption, for even $L \geq 4$, there is little margin to improve Lemma~\ref{le:uncertainty}: at best, one can hope to replace the strict inequality in~\eqref{eq:sparsityassumption0} with a large one. Can this actually be done ? This is partially answered by the following results:
\begin{lemma}
\label{le:sharpuncertaintydisjstrict} Assume that $\tA \not\equiv \A$, that~\eqref{eq:ksparsefilters} holds with \begin{equation}
\label{eq:sparsityassumption1}
2 k \cdot \Delta(\tA,\A) = L
\end{equation}
and that the filters in $\A$ have {\em disjoint supports}~\eqref{eq:disjointsupports}.
Then, either $\lo{\tA} > \lo{\A}$, or each row of $\tA$ is obtained by permuting pairs of distinct filters $\a_{ij}$, $\a_{ij'}$ from the corresponding row of $\A$ such that $\a_{ij}-\a_{ij'}$ is proportional to a modulated and translated Dirac comb with $2k$ spikes. 
\end{lemma}

For filter matrices with a single row, since $\tA \not\equiv\A$ means that the filters $\ta_{1j}$ are permuted versions of $\a_{1j}$, we obtain
\begin{coro}
\label{cor:singlerow} Consider $\A$ with a single row ($M=1$). Assume that $\tA \not\equiv \A$, that~\eqref{eq:ksparsefilters} holds with \begin{equation}
\label{eq:sparsityassumption2}
2 k \cdot \Delta(\tA,\A) = L
\end{equation}
and that the filters in $\A$ have {\em disjoint supports}~\eqref{eq:disjointsupports}. Then $\lo{\tA} > \lo{\A}$. 
\end{coro}

\subsection{Excessive pessimism?}
The counter-examples built in Lemmata~\ref{le:sharpuncertainty}-\ref{le:sharpuncertaintydisj}, which are associated to Dirac combs, are highly structured. They provide worst case well-posedness bounds, but existing probabilistic versions of uncertainty principles (see, e.g., the nice survey~\cite{springerlink:10.1007/s00041-008-9042-0}) lead us to conjecture that if the sparse filters in $\A$ are drawn at random (e.g. from Bernoulli-Gaussian distribution), the uniqueness guarantee of Theorem~\ref{th:main} will hold except with small probability $O(L^{-\beta})$, provided that $k < c(\beta) L/\log L$, for large $L$. 
This is 
 left to further theoretical investigation. 

%
%
%
\section{Numerical experiments}
\label{sec:algorithm}

The results achieved so far are theoretical well-posedness guarantee, but do not quite provide algorithms to compute the potentially unique (up to global permutation) solution of the frequency permutation problem. We conclude this paper with the description of a relatively naive optimization algorithm, an empirical assessment of its performance with Monte-Carlo simulations, and a discussion of how this compares with the theoretical uniqueness guarantees achieved above.
 
\subsection{Proposed combinatorial algorithm}
Given a ``permuted'' matrix $\tA$, one wishes to find a set of frequency permutations yielding a new matrix $\wA$ with minimum $\lp$ norm. 

The proposed algorithm starts from $\wA_{0} = \tA$. Given $\wA_{n}$,  a candidate matrix $\wA_{n+1,\pi}$ can be obtained by applying a permutation $\pi$ at frequency $\omega_{n} \equiv n\ [\texttt{mod}\ L]$. Testing each possible permutation $\pi$ and retaining the one $\pi_{n}$ which minimises $\|\wA_{n+1,\pi}\|_{p}$ yields the next iterate $\wA_{n+1} := \wA_{n+1,\pi_{n}}$. The procedure is repeated until the $\lp$ norm $\wA_{n}$ ceases to change. Since there is a finite number of permutations to try, the stopping criterion is met after sufficiently many iterations.

\subsection{Choice of the $\lp$ criterion}
In theory, it could happen that the stopping criterion is only met after a combinatorially large number of iterations. However, the algorithm stops much sooner in practice. In fact, if we were to use the $\ell^{0}$ norm, the algorithm would typically stop after just one iteration, because the $\ell^{0}$ norm attains its maximum value $M \times N \times L$ for most frequency permutations except a few very special ones. For this reason, we chose to test the algorithm using $\lp$ norms $p>0$, which are not as ``locally constant'' as the $\lzero$ norm. To our surprise, the experiments below will show that the best performance is not achieved for small $p$, but rather for $p=2-\epsilon$ with small $\epsilon>0$. For $p=0$ and for $p\geq 2$, the algorithm indeed completely fails.

\subsection{Monte-Carlo simulations}
For various filter length $\filterlength$, sparsity levels $k$ and dimensions $\nrofchn$, $\nrofsrc$, random sparse filter matrices $\A$ made of independent random $k$-sparse filters were generated. Each filter was drawn by choosing: a) a support of size $k$ uniformly at random; b) i.i.d. Gaussian coefficients on this support. 

For each configuration $(L,k,\nrofchn,\nrofsrc)$, 200 such random matrices $\A$ were drawn. For each $\A$, independent random frequency permutations were applied to obtain  $\tA$. The algorithm was then applied to obtain $\wA$. The performance was measured using the SNR between $\A$ and the best permutation of $\wA$. 

\begin{figure}[t]
\centering
\includegraphics[width=\textwidth*4/5,height=\textwidth/2]{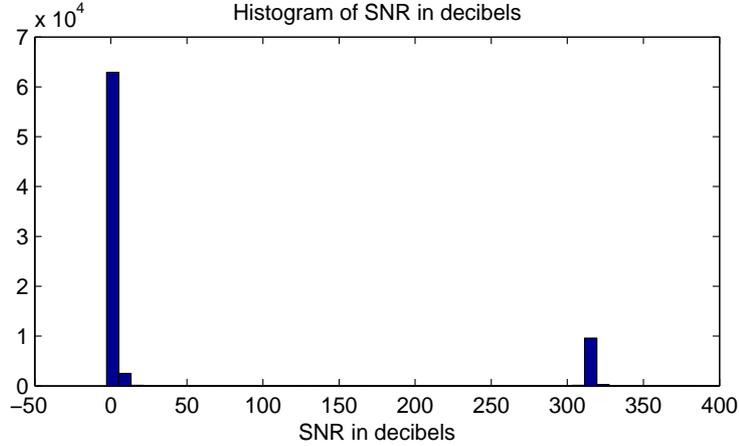}
\caption{Histogram of SNR between best permutation of $\wA$ and original $\A$}
\label{fig:SNRHist}
\end{figure}

Figure~\ref{fig:SNRHist} shows the histogram of SNR values achieved for $L=31$, $1 \leq k \leq L$, $M \in \{1,2\}$, $N\in\{2,3,4\}$, $p=1$. 
It shows that the algorithm either completely succeeds up to machine precision (SNR above 300 dB) or completely fails (SNR of the order of 0 dB). For this reason, in the rest of the experiments the estimation was considered a success when the SNR exceeded 100 dB.

\subsection{Role of the $\ell^{p}$ norm}

\begin{figure}[t]
\centering
\includegraphics[width=\textwidth*4/5,height=\textwidth/2]{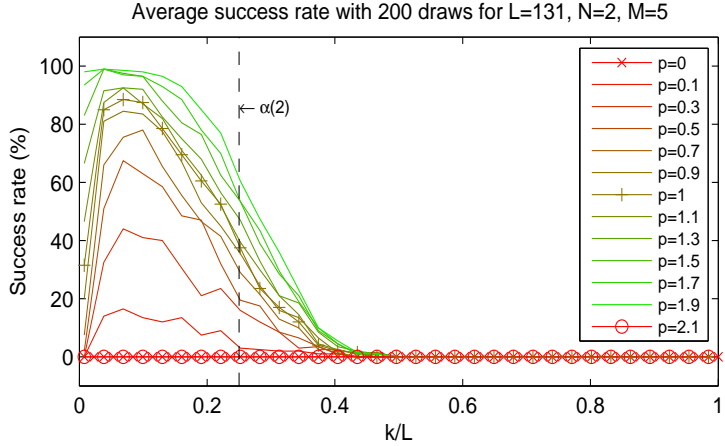}
\caption{Filter recovery success as a function of $p$, $0\leq p\leq 1.9$}
\label{fig:roleofp}
\end{figure}

Figure~\ref{fig:roleofp} displays the success rate as a function of the relative sparsity $k/L$, for various choices of the $\lp$ criterion,  with filters of prime length $L=131$, $N=2$ sources and $M=5$ channels. The vertical dashed line indicates the threshold $k/L \leq \alpha(2)$ associated with the well-posedness guarantee (using an $\ell^0$ criterion) of Theorem~\ref{th:main}.
Surprisingly, one can observe that the success rate increases when $0<p<2$ is increased. The maximum success rate is achieved when $p=2-\epsilon$ with small $\epsilon>0$. 

Beyond the well-posedness regime suggested by the theory (i.e., to the right of the vertical dashed line) the algorithm can succeed, but at a rate that rapidly decreases when the relative sparsity $k/L$ increases. In the regime where the problem is proved to be well-posed, the proposed algorithm is often successful but can still fail to perfectly recover the filters, especially --and surprisingly-- for small values of $k$. This phenomenon is strongly marked for $p<1$ and essentially disappears for $p>1$. It remains an open question to determine the respective roles of the $\lp$ criterion and of the naive greedy optimization algorithm in this limited performance when the problem is well-posed with respect to the $\lzero$ norm.

\subsection{Role of the filter length $L$}
\begin{figure}[t]
\centering
\includegraphics[width=\textwidth*4/5,height=\textwidth/2]{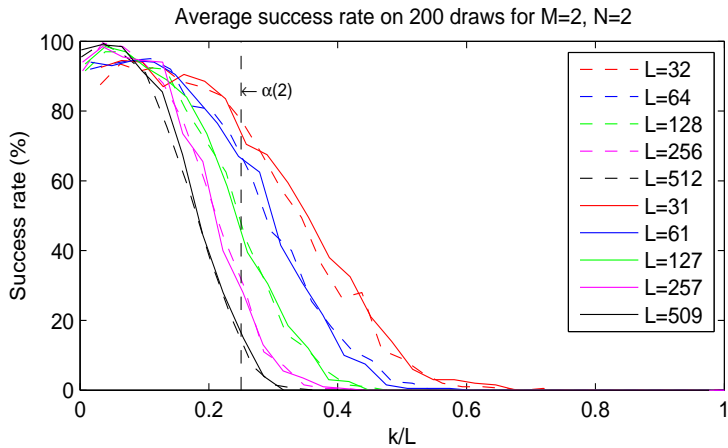}
\caption{Filter recovery success as a function of $L$, for $p=1.9$}
\label{fig:roleofL}
\end{figure}

Figure~\ref{fig:roleofL} shows the results for different $L$ values 
with $p=1.9$, $M=N=2$. 
One can see that the average performance does not seem to depend on whether $L$ is prime or not.
As $L$ increases, the performance for ``small'' $k/L$ slightly increases, but the success rate degrades for ``large'' $k/L$ close to $\alpha(2)$. 

\subsection{Role of the number of channels $M$}
\begin{figure}[t]
\centering
\includegraphics[width=\textwidth*4/5,height=\textwidth/2]{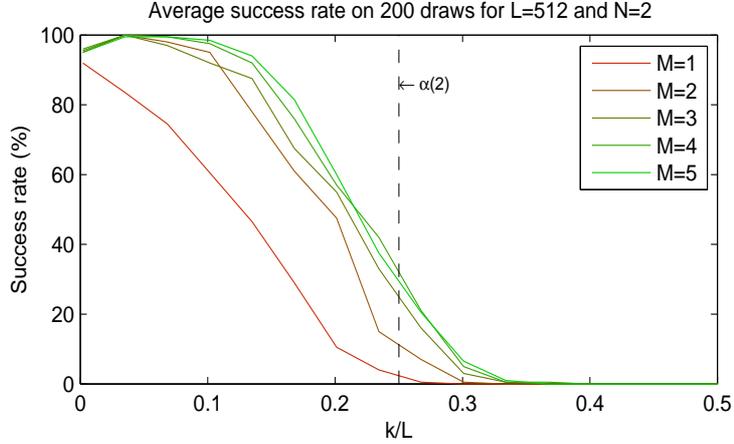}
\caption{Filter recovery success as a function of $M$, for $p=1.9$}
\label{fig:roleofM}
\end{figure}
Figure~\ref{fig:roleofM} shows the results for increasing numbers of channels $M$, with a filter length $L=512$, $N=2$ sources, $p=1.9$. 
One can observe that the success rate substantially increases when $M$ is increased from $M=1$ to $M=2$, and slightly increases as $M$ further increases. 
Although the worst-case well-posedness guarantees are the same, the algorithm seems to benefit from added filter diversity across channels.

\subsection{Role of the number of sources $N$}
\begin{figure}[t]
\centering
\includegraphics[width=\textwidth*4/5,height=\textwidth/2]{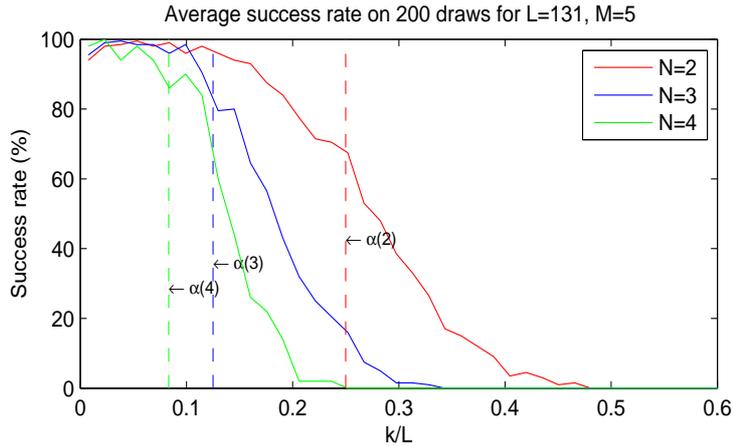}
\caption{Filter recovery success as a function of $N$, for $p=1.9$}
\label{fig:SuccessRate}
\end{figure}
Figure~\ref{fig:SuccessRate} shows the success rate as a function of the relative sparsity $k/L$,  for $\nrofsrc\in \{2,3,4\}$, with $L=31$, $\nrofchn =5$ with $p=1.9$. The well-posedness limits $k/L \leq \alpha(N)$ associated to Theorem~\ref{th:main} are indicated with vertical dashed lines.
The empirical curves confirm that the algorithm can still succeeed beyond the worst-case well-posedness guarantees, but with a rapidly decreasing rate of success. When the well-posedness guarantees hold, the algorithm can fail, but its rate of success is high when the relative sparsity is sufficiently small compared to the bound provided by Theorem~\ref{th:main}.

\subsection{Computation time}
The algorithm evaluates the $\ell^p$ norm of the $N!$ permutations of the sources for each of the $L$ frequencies.

To evaluate the $\ell^p$ norm of the filters, the permuted frequency coefficients have to be transformed back into the time domain by inverse Discrete Fourier transform. For each filter, the cost of the Discrete Fourier Transform through a Fast Fourier Transform is $\mathcal{O}(L \log_{2} L)$. There are $MN$ filters and hence the cost of $\ell^p$ norm evaluation for a given configuration of sub-bands is $\mathcal{O}(MNL\log_{2} L)$. 

Hence, the complexity of each sweep through the set of all frequencies is $\mathcal{O}(N!MNL^2 \log_{2} L)$. This is rather expensive because the computational cost grows in factorial with the number of sources and in square with the filter length, but it is tractable for small problem sizes and very efficient compared to the brute force approach that would require $\mathcal{O}( (N!)^{L-1} MNL \log_{2} L)$ operations to test all $(N!)^{L-1}$ possible permutations up to a global permutation.
 
\begin{figure}[t]
\centering
\includegraphics[width=\textwidth*4/5,height=\textwidth/2]{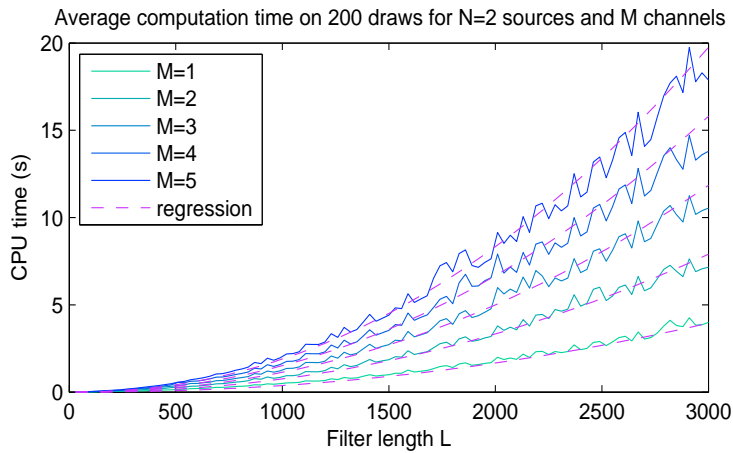}
\caption{Computation time of the permutation solving algorithm depending on the length $L$ of the filter}
\label{fig:cpu}
\end{figure}

Figure \ref{fig:cpu} shows the average computation time over $200$ trials for various filter length. The red dashed line corresponds to its prediction using the theoretical cost estimation as $C \times L^2 \log_{2} L$ with $C \approx 40$ nanoseconds.


%
%


\section{Conclusions}

It is now well known that a sufficient sparsity assumption can be used to make under-determined linear inverse problems well-posed: without the sparsity assumption, the problem admits an affine set of solutions, which intersects at only one point with the set of sparse vectors. Besides this well-posedness property, a key factor that has lead to the large deployment of sparse models and methods in various fields of science is the fact that a convex relaxation of the NP-hard $\ell^{0}$ minimization problem can be guaranteed to find this unique solution under certain sparsity assumptions. The availability of efficient convex solvers then really makes the problem tractable.
 
 The problem considered in this paper is not a linear inverse problem. Even though it is a simplification of the original permutation {\em and scaling} problem arising from signal processing, it remains a priori a much harder problem than linear inverse problems in terms of the structure of the solution set: each solution comes with a herd of solutions that are equivalent up to a global permutation. 
 
 The fact that we managed to obtain well-posedness results in this context is encouraging, but this is at best the beginning of the story: even if the solution is unique, how do we efficiently compute it? Can one hope to extend these results to the original permutation and scaling problem? Why does the proposed naive algorithm perform better for $p>1$? Answers to these questions are likely to have an impact in fields such as blind source separation with sparse multipath channels.

\appendices


%

\section{Proof of Theorem 1}

First, notice that for each frequency $\omega$ and channel $i$, permutations preserve the equality $\sum_{j} a_{ij}[\omega]=\sum_{j} \tilde {a}_{ij}[\omega]$. Thus, the same holds in the time-domain : 
$\sum_{j} a_{ij}=\sum_{j} \tilde {a}_{ij}$.
By the disjoint supports hypothesis and the quasi-triangle inequality for $\lp$ quasi-norms we have
\begin{equation}
\sum_{j} \|a_{ij}\|_p^p
=
\|\sum_{j} a_{ij}\|_p^p
=\|\sum_{j} \tilde a_{ij}\|_p^p\leq\sum_{j}  \|\tilde a_{ij}\|_p^p.
\end{equation}
We conclude by summing over all channels $i$.
%
%


\section{Proof of Lemma~\ref{le:BiStochastic} and Corollary~\ref{cor:perms}}
\label{app:PfMain}

Let us consider the matching count matrix $\mathbf{C}$ with entries 
\[
C_{jn} := \sharp \{0 \leq \ell < L: \sigma_{\ell}(j)=n\},\quad 1 \leq j,n\leq N.
\]
Since $\sum_{j} C_{jn} = \sum_{n} C_{jn} = L$ we have $\mathbf{C} = L \cdot \mathbf{B}$ where $\mathbf{B}$ is bi-stochastic. 

A weakened version of Lemma~\ref{le:BiStochastic}, with $2\alpha''(N) = \frac 1 {1+(N-1)^2}$, can be obtained by combining the Birckhoff - Von Neumann theorem and Carath\'eodory theorem. 
 \begin{theorem}[Birkhoff - Von Neumann  Theorem, \cite{birkhoff1946three,von1953certain}]
Every bi-stochastic matrix is in the convex hull of permutation matrices. 
 \end{theorem}
 \begin{theorem}[Carath\'eodory Theorem \cite{berger1977geometrie}]
 Let $X$ be a non empty subset of an affine space of dimension $n\geq 1$. Then every element of the convex hull of X is a convex combination of $p$ elements of $X$, with $p\leq n+1$.
 \end{theorem}
 
The set of bi-stochastic matrices is an affine subspace of $\mathbb{R}^{N^{2}}$. It is defined by $2N$ equations, but these equations are linearly dependent since the sum of the sum over rows is the sum of the sum over columns. Hence its affine dimension is $n \leq N^{2}-(2N-1) = (N-1)^{2}$, and we conclude from Carath\'eodory's theorem that every bi-stochastic matrix is a convex combination of at most $(N-1)^{2}+1$ permutation matrices. One of the coefficients of this combination must therefore exceed $1/(1+(N-1)^{2})$, and this leads to a version of Lemma~\ref{le:BiStochastic} with $2\alpha''(N)=\frac 1 {1+(N-1)^2}$, as claimed.

Yet, this bound is suboptimal. The optimal bound in Lemma~\ref{le:BiStochastic} follows from Hall's Marriage Theorem, which by the way is also a key ingredient in the proof of the Birkhoff-Von Neumann theorem.
 \begin{theorem}[Hall's Marriage Theorem \cite{Hall:1935aa,oxleymatroid} ]
Let $(A_j)_{j\in J}$ be a family of subsets of a set finite $S$. There exists a bijection $\pi:J\rightarrow S$  such that $\pi(j)\in A_{j}$ for all $J$ if, and only if, for all $E \subset J$ $$\sharp \cup_{j\in E} A_j\geq \sharp E$$
\end{theorem}
The bijection $\pi$ is often referred to as a \textit{transversal} for $S$. 

\begin{proof}[Proof (Lemma~\ref{le:BiStochastic})]
For shortness of notation we write $\alpha$ for $\alpha(N)$. Define the sets $J = S = \llbracket 1,\ N\rrbracket$, and $A_{j} := \{n: B_{jn}\geq 2\alpha\}$, $j \in J$, and consider the property \[
\mathcal{P}_{k}: \forall E \subset J, \sharp E \leq k \Rightarrow \sharp \cup_{j \in E} A_{j} \geq \sharp E.
\]
We wish to prove that $\mathcal{P}_{k}$ holds true for all $1 \leq k \leq N$: then, by Hall's Marriage Theorem, there exists a bijection $\pi: j \to \pi(j)$ such that $\pi(j) \in A_{j}$ for all $j$, yielding in turn the permutation matrix $\mathbf{P}$ with ones at the entries $(j,\pi(j))$. We proceed by contradiction: assume that $\mathcal{P}_{N}$ does not hold true. Since $\mathcal{P}_{1}$ holds true, without loss of generality, for some $1 \leq k_{0} < N$: 
\[
\sharp \cup_{1 \leq j \leq k_{0}} A_{j} \geq k_{0},\ \mbox{and}\ 
\sharp \cup_{1 \leq k \leq k_{0}+1} A_{j} \leq k_{0}.
\] 
Hence, without loss of generality: 
\[
\cup_{1 \leq k \leq k_{0}} A_{j} = \llbracket 1\  k_{0}\rrbracket \supset A_{k_{0}+1}.
\]
It follows that for $n>k_{0}$ and $j \leq k_{0}+1$, we have $n \notin A_{j}$, hence $B_{jn} < 2\alpha$. Now we use the bi-stochasticity of $\mathbf{B}$ ($\sum_{j} B_{jn} = \sum_{n} B_{jn}=1$, $B_{jn} \geq 0$) to obtain
\begin{eqnarray*}
k_{0} & \geq & \sum_{n \leq k_{0}} \sum_{j \leq k_{0}+1} B_{jn}
 =  \sum_{j \leq k_{0}+1} \sum_{n \leq k_{0}} B_{jn}\\
& = & \sum_{j \leq k_{0}+1} \Big(1-\sum_{n>k_{0}} B_{jn}\Big)
 >  \sum_{j \leq k_{0}+1} \left(1-(N-k_{0})2\alpha \right)\\
 & = & (k_{0}+1) (1-(N-k_{0})2\alpha)\\ 
 &=& k_{0}+\Big(1-(k_{0}+1)(N-k_{0})2\alpha\Big).
\end{eqnarray*}
This implies $2\alpha > 1/(k_{0}+1)(N-k_{0})$. However, this yields a contradiction, since a simple functional study shows that 
\[
\max_{1 \leq k_{0}<N} \frac{1}{(k_{0}+1)(N-k_{0})} = 2\alpha.
\]
 \end{proof}\qed 

Equipped with Lemma~\ref{le:BiStochastic}, we can now prove Corollary~\ref{cor:perms}.
\begin{proof}[Proof (Corollary~\ref{cor:perms})]
Since $\mathbf{C} = L \cdot \mathbf{B}$ where $\mathbf{B}$ is bi-stochastic, there is a permutation $\pi$ such that $C_{j\pi(j)} \geq 2L\alpha(N)$.
 \end{proof}\qed

We conclude this section by showing the sharpness of Corollary~\ref{cor:perms} through the construction of permutations that reach the bound.
Consider $N$ an integer, and $k_{0} := N/2$ ($N$ even) or $k_{0} := (N-1)/2$ ($N$ odd). Let $L$ be a multiple of $(k_{0}+1)(N-k_{0})$. 
Consider the $L \times N$ matrix:
\[
\mathbf{\Sigma} :=
\left[
\begin{array}{ccccc|ccc}
\mathbf{1} & \mathbf{U} & \mathbf{k_{0}} & \ldots  & \mathbf{2} & \times & \ldots & \times\\
\mathbf{2} & \mathbf{1} & \mathbf{U} & \ddots & \ddots & \vdots & & \vdots \\
\vdots & \ddots & \ddots & \ddots & \ddots & \vdots & & \vdots \\
\mathbf{k_{0}} & \ddots & \ddots & \ddots & \mathbf{U} & \vdots & & \vdots \\
\mathbf{U} & \mathbf{k_{0}} & \ldots & \mathbf{2} & \mathbf{1} & \times & \ldots & \times
\end{array}
\right]
\]
where: a) the left part, of size $L \times (k_{0}+1)$, is filled with the column vectors $\mathbf{i} \in \R^{L/(k_{0}+1)}$ made of constant entries equal to the integer $1 \leq i \leq k_{0}$ and the vector $\mathbf{U} \in \R^{L/(k_{0}+1)}$ made of the vertical concatenation of the $N-k_{0}$ column vectors $\mathbf{j} \in \R^{L/(k_{0}+1)(N-k_{0})}$ with constant entries $k_{0}+1 \leq j \leq N$; b) the rows of the the right part, of size $L \times (N-k_{0}-1)$, include exactly once each integer $1 \leq \ell \leq N$ which does not already appear in the corresponding row of the left part.
By construction, the $L$ rows of the matrix $\mathbf{\Sigma}$ are associated to $L$ permutations $\sigma_{\ell}$. 
We now show that, for any global permutation $\pi$, there is at least one column $1 \leq j \leq k_{0}+1$ such that 
\[
\sharp \{\ell: \sigma_{\ell}(j) = \pi(j)\} \leq L/(k_{0}+1)(N-k_{0}) = L 2\alpha(N).
\]
Applying again the pigeonhole principle yields: among the $k_{0}+1$ indices $j$ to consider, at least one, $j^{\star}$, must be mapped to an integer $\pi(j^{\star}) \geq k_{0}+1$. By construction, the columns of $\Sigma$ are such that column $j^{\star}$ contains at most (in fact: exactly) $L/(k_{0}+1)(N-k_{0})$ instances of the value $\pi(j^{\star})$. 

\section{Proof of  Theorem~\ref{th:main}}
We can now conclude the proof of Theorem~\ref{th:main}. By Corollary~\ref{cor:perms}, there is a permutation $\pi$ such that for each $j$, we have $\|\mathbf{F}(\ta_{ij}-\a_{i\pi(j)})\|_{0} \leq L(1-2\alpha(N))$, hence $\Delta(\tA,\A|\pi) \leq L(1-2\alpha(N))$ and finally $\Delta(\tA,\A) \leq L(1-2\alpha(N))$. Combined with the assumption $k \leq L \alpha(N)$, we obtain $2k+\Delta \leq L$, and we conclude thanks to Lemma~\ref{le:uncertainty}.

\section{Proof of Lemmata~\ref{le:uncertainty}~\ref{le:sharpuncertainty}~\ref{le:sharpuncertaintydisj}~\ref{le:sharpuncertaintydisjstrict}}
%
%
We prove Lemma~\ref{le:uncertainty} first, then the statements of Lemma~\ref{le:sharpuncertainty} in the following order: 1), 3), 2). We begin by some notations and fact regarding Dirac combs.

\subsection{Dirac combs}
\label{sec:diraccombs}
Let $p,q \geq 1$ be two integers and $L=pq$ their product. The unit Dirac comb with $q$ spikes and of step $p$, denoted  $\dirac_{p}$, is the vector of $\C^{L}$ defined by $\dirac_{q}[t]=1/\sqrt{q}$ if $t \equiv 0 [p]$, $\dirac_{q}[t]=0$ otherwise. Its Fourier transform is the unit Dirac comb with $p$ spikes and of step $q$: $\F\dirac_{q} = \dirac_{p}$.
For $0\leq n<p$ an integer translation index and $0 \leq m<q$ an integer modulation index, one can define
the translated and modulated Dirac comb $\dirac_{q,n,m} = T_{n} M_{m} \dirac_{q}$ where $T_{n}$ is the circular shift by $n$ samples, and $M_{m}$ is the frequency modulation $(M_{m}u)[t] := u[t] \cdot e^{2i\pi mt/L}$.  One can check that the collection $\{\dirac_{q,n,m}\}_{0\leq n<p,0\leq m<q}$ is an orthonormal basis of $\C^{L}$. 

\subsection{Proof of Lemma~\ref{le:uncertainty}}
Let $\pi_0$ be the permutation such that $\Delta_{0}(\A,\tA) 
= \min_{\pi\in\mathfrak S_N} \Delta_{0}(\A,\tA|\pi)$. By abuse of notation we still denote $\A$ the matrix obtained by applying $\pi_{0}$ to permute the columns of the original filter matrix. 
For each channel $i$ and a source index $j$ such that $\a_{ij}=\ta_{ij}$ we obviously have $\lo{\a_{ij}} \leq \lo{\ta_{ij}}$. Now, since $\Delta_{0} \geq 1$ we have $\tA \not\equiv \A$ hence there exists a pair $i,j$ such that $\ta_{ij} \neq \a_{ij}$. By the $\lzero$ Dirac-Fourier uncertainty principle  \cite[Theorem 1]{elad2002generalized}, for any vector $u \in \C^{\filterlength}$ we have $\lo{u} \lo{\F u} \geq \filterlength$. Hence, by the hypothesis $k < \filterlength/(2\Delta_{0})$ we have
\begin{eqnarray}
\lo{\a_{ij}} + \lo{\ta_{ij}}
& \geq &
\lo{\ta_{ij}-\a_{ij}}\label{ineq:1}\\
& \geq & \filterlength/\lo{\F(\ta_{ij}-\a_{ij})}\label{ineq:2}\\
& \geq & \filterlength/\Delta_{0} > 2k\label{ineq:3}\\
& \geq & \lo{\a_{ij}} +\lo{\a_{ij'}}\label{ineq:4}
\end{eqnarray}
where $j'$ is an arbitrary source index.
Hence for every $i,j$ such that $\ta_{ij} \neq \a_{ij}$ and any $j'$, $\lo{\ta_{ij}} > \lo{\a_{ij'}}$, and we obtain
\[
\lo{\ta_{ij}}> \max_{j'} \lo{\a_{ij'}}  \geq \lo{\a_{ij}}.
\]
Overall, we have shown that $\lo{\tA} > \lo{\A}$.

When $L$ is prime, a stronger uncertainty principle $\lo{u}+\lo{\mathbf{F}u}\geq L+1$ holds \cite{Tao:2005aa}. Hence, under the assumption $2k + \Delta_{0} \leq L$ we can replace~\eqref{ineq:2}-\eqref{ineq:3} with 
\[
\ldots \geq L+1-\lo{\F (\ta_{ij}-\a_{ij})} \geq L+1-\Delta_{0} > 2k
\]
to reach the same conclusion. 


%


\subsection{Proof of Lemma~\ref{le:sharpuncertainty}}
We shall simply build an example where $\A = [\alpha,\ \beta]$ is a $1 \times 2$ matrix of filters. Extensions to $\A$ an $\nrofchn \times \nrofsrc$ matrix are trivial by adding mutually distinct sparse columns that are distinct from $\alpha$ and $\beta$, and duplicating the first row. 

We exploit Dirac combs as described in Appendix~\ref{sec:diraccombs}.
Define $a = \dirac_{k,0,0}$, $b = -\dirac_{k,L/2k,0}$. The filters $a$ and $b$ have disjoint support and satisfy $\lo{a} = \lo{b} = k$. Since $a-b= \sqrt2\ \dirac_{2k,0,0}$ 
we have $a[\omega] = b[\omega]$ whenever $\omega \not\equiv 0 [2k]$. Hence, permuting the Fourier transforms of $a$ and $b$ on the $L/2k$ frequencies $\{ \omega = 2kr, 0 \leq r < L/2k \}$ yields $\tilde{a} = b$ and $\tilde{b} = a$.
Given any $u \in \C^L$ we define perturbations $\alpha$ and $\beta$ of $a$ and $b$
\[
\begin{cases}
\alpha&:=a+u\\
\beta&:=b+T_{L/2}  u
\end{cases}
\]
with $T_{L/2}$  a circular shift. Noticing that for $\omega=2kr$ 
\[
(T_{L/2}u) [\omega]=e^{\frac{2i\pi (L/2)\omega } {L}} u[\omega]=e^{2i\pi kr} u[\omega] = u[\omega]
\]
we obtain that, after permuting the Fourier transforms of $\alpha$ and $\beta$ at the frequencies $\omega = 2kr$, $0 \leq r < L/2k$,
\[
\begin{cases}
\tilde \alpha&=b+u\\
\tilde \beta&=a+T_{L/2} u
\end{cases}
\]
We choose the vector $u$ to be zero everywhere with two exceptions 
$u[0]:=-a[0]$, $u[\textstyle \frac L{2k}]:=-b[\textstyle \frac L{2k}]$. 
Since $T_{L/2}u \neq u$ and $a \neq b$, we have $\{\alpha,\beta\}\neq \{\tilde\alpha,\tilde \beta\}$ and $\tA \not\equiv\A$. 
Moreover, $\Delta_{0}(\tA,\A)=\Delta_{1}(\tA,\A)=L/2k$.

Lastly, all considered vectors have $k$ entries of equal magnitude, hence $\|\alpha\|_0=\|\beta\|_0=\|\tilde\alpha\|_0=\|\tilde\beta\|_0=k$, and for any $0<p \leq \infty$ $\|\alpha\|_p=\|\beta\|_p=\|\tilde\alpha\|_p=\|\tilde\beta\|_p$. In particular, $\|\tA\|_{p} = \|\A\|_{p}$, $0 \leq p \leq \infty$. 

\subsection{Proof of Lemma~\ref{le:sharpuncertaintydisjstrict}}
We repeat the construction of the proof of Lemma~\ref{le:sharpuncertainty} starting from the Dirac combs $a=\dirac_{k',0,0}$, $b=-\dirac_{k',L/2k',0}$. Since $k'< k \leq L/2$, we have $\ell:= k-k' \leq L/2-k'$ hence we can choose an $\ell$-sparse vector $u$ which support is outside the support of  $\dirac_{2k'}$ and such that $T_{L/2}u$ and $u$ have disjoint supports. The four vectors $\{a,b,u,T_{L/2}  u\}$ have mutually disjoint supports, hence $\alpha$ and $\beta$ have disjoint supports, $\{\alpha,\beta\}\neq \{\tilde\alpha,\tilde \beta\}$ and $\tA \not\equiv\A$. Moreover, $\Delta_{0}(\tA,\A) = \Delta_{1}(\tA,\A) = L/2k$. Lastly, we have $\lo{\alpha} = \lo{\beta} = \lo{\tilde{\alpha}} = \lo{\tilde{\beta}} = k'+\ell = k$, and the $\lp$ norms of these vectors are also equal, hence $\|\tA\|_{p} = \|\A\|_{p}$, $0 \leq p \leq \infty$. 

\subsection{Proof of Lemma~\ref{le:sharpuncertaintydisj}}
As in the proof of Lemma~\ref{le:uncertainty} we consider $\A$ the permuted matrix associated to the optimal permutation $\pi_{0}$. Using the inequality $2k \leq L/\Delta_{1} \leq L/2\Delta_{0}$ instead of $2k < L/\Delta_{0}$ we repeat the steps~\eqref{ineq:1}-\eqref{ineq:4} to obtain $\lo{\ta_{ij}} \geq \lo{\a_{ij'}}$ for any $j \in E_{i} := \{j, \a_{ij} \neq \ta_{ij}\}$ and any $j'$. As a result $\lo{\ta_{ij}} \geq \lo{\a_{ij}}$ for all $i,j$.
The assumption that $\lo{\tA} = \lo{\A}$ 
implies that $\lo{\ta_{ij}} = \lo{\a_{ij}}$ for all $i,j$.

By assumption, $\tA \not\equiv \A$ hence there are indices $i,j$ such that $\a_{ij} \neq \ta_{ij}$.
For such $i,j$, since $\lo{\a_{ij}} = \lo{\ta_{ij}}$, each inequality in~\eqref{ineq:1}-\eqref{ineq:4} (the inequality $L/\Delta_{0}>2k$ being replaced with $L/\Delta_{1}\geq 2k$) must be indeed an equality. This implies that: $\lo{\a_{ij}} = \lo{\ta_{ij}}=k$; $2k$ divides $L$ and $\Delta_{1}=L/2k$; the nonzero vector $b_{ij} := \ta_{ij}-\a_{ij}$ must be an equality case of the $\lzero$ uncertainty principle with $\lo{b_{ij}}=2k$ and $\lo{\mathbf{F}b_{ij}}=L/2k$. As a result \cite{Tao:2005aa} $b_{ij}$ is a scaled, modulated and translated version of the Dirac comb $\dirac_{2k}$ made of $2k$ Diracs  spaced every $L/2k$ samples: there exists a scalar $\gamma_{ij} \neq 0$, and two integers $0 \leq n_{ij} < L/2k$, $0 \leq m_{ij} < 2k$ such that
\[
b_{ij} = \gamma_{ij} \cdot \dirac_{2k,n_{ij},m_{ij}}.
\]
Moreover since $\lo{\a_{ij}} = \lo{\ta_{ij}} = k$ and
$\|\ta_{ij}-\a_{ij}\|_0=2k$, the filters $\ta_{ij}$ and $\a_{ij}$ have disjoint supports of size $k$. Hence, they are the restriction of $b_{ij}$ (resp. of $-b_{ij}$) to their respective supports. 

Now, define 
\[
E_{i,n,m} := \{j \in E_{i}, n_{ij}=n,m_{ij} =m\}.
\]
As observed in the proof of Theorem~1, the equality $\sum_{j} \a_{ij} = \sum_{j} \ta_{ij}$ holds, implying $\sum_{j \in E_{i}} b_{ij} = \sum_{j} b_{ij} = 0$. Taking inner products with the Dirac comb orthonormal basis $\dirac_{2k,n,m}$, $0 \leq n<L/2k$, $0\leq m<2k$, yields 
\begin{eqnarray}
\sum_{j \in E_{i,n,m}} \gamma_{ij} &=& 0,\label{eq:sumgamma}
\end{eqnarray}
Since $\gamma_{ij} \neq 0$, whenever $E_{i,n,m}$ is not empty it contains at least two distinct indices. 

By the {\em disjoint support assumption}:
for  $j,j' \in E_{i,n,m}$, $j\neq j'$, the original filters $\a_{ij}$ and $\a_{ij'}$ have disjoint supports. Moreover, we know that these supports are subsets of the support of $\dirac_{2k,n,m}$ which is of size $2k$, hence
\[
 \sharp E_{i,n,m} \cdot k = \lo{\sum_{j \in E_{i,n,m}} \a_{ij}} \leq 2k.
\]
Hence, whenever $E_{i,n,m}$ is not empty, it contains {\em exactly} two distinct elements: $E_{i,n,m} = \{j,j'\}$ where $j \neq j'$. 

Further, observe that: a) $\a_{ij}$ and $\a_{ij^\prime}$  have disjoint supports of size $k$ which are subsets of the support of size $2k$ of $\dirac_{2k,n,m}$; b) $\a_{ij}$ and $\ta_{ij}$ have the same property. As a result, $\ta_{ij}$ and $\a_{ij'}$ have the same support, which is disjoint from that of $\a_{ij}$. Similarly, $\a_{ij}$ has the same support as $\ta_{ij'}$. Finally, Eq.~\eqref{eq:sumgamma} can be rewritten $\gamma_{ij} + \gamma_{ij'}=0$, and implies $b_{ij} + b_{ij'} = 0$, that is to say $\ta_{ij} +\ta_{ij'} = \a_{ij'} + \a_{ij}$. We conclude that $\ta_{ij} = \a_{ij'}$ and $\ta_{ij'}=\a_{ij}$.

\section*{Acknowledgment}
This work was supported by the EU Framework 7 FET-Open project
FP7-ICT-225913-SMALL: Sparse Models, Algorithms and Learning for
Large-Scale data, and by Agence Nationale de la Recherche (ANR), project ECHANGE (ANR-08- EMER-006).


\bibliographystyle{unsrt}
\bibliography{bigbib}

\begin{thebibliography}{10}

\bibitem{5352256}
C.R. Berger, Shengli Zhou, J.C. Preisig, and P.~Willett.
\newblock Sparse channel estimation for multicarrier underwater acoustic
  communication: From subspace methods to compressed sensing.
\newblock {\em Signal Processing, IEEE Transactions on}, 58(3):1708 --1721,
  March 2010.

\bibitem{4518398}
M.~Sharp and A.~Scaglione.
\newblock Application of sparse signal recovery to pilot-assisted channel
  estimation.
\newblock In {\em Acoustics, Speech and Signal Processing, 2008. ICASSP 2008.
  IEEE International Conference on}, pages 3469 --3472, April 2008.

\bibitem{5454399}
W.U. Bajwa, J.~Haupt, A.M. Sayeed, and R.~Nowak.
\newblock Compressed channel sensing: A new approach to estimating sparse
  multipath channels.
\newblock {\em Proceedings of the IEEE}, 98(6):1058 --1076, June 2010.

\bibitem{benestybook}
Jacob Benesty, M.~Mohan Sondhi, and Yiteng~(Arden) Huang.
\newblock {\em Springer Handbook of Speech Processing}.
\newblock Springer-Verlag New York, Inc., Secaucus, NJ, USA, 2007.

\bibitem{934137}
Binning Chen and A.P. Petropulu.
\newblock Frequency domain blind mimo system identification based on second and
  higher order statistics.
\newblock {\em Signal Processing, IEEE Transactions on}, 49(8):1677 --1688, Aug
  2001.

\bibitem{2010_AcademicPress_BookBSS}
Pierre Comon and Christian Jutten, editors.
\newblock {\em Handbook of Blind Source Separation, Independent Component
  Analysis and Applications}.
\newblock Academic Press, 2010.

\bibitem{97}
M.S. Pedersen, J.~Larsen, U.~Kjems, and L.C. Parra.
\newblock A survey of convolutive blind source separation methods.
\newblock {\em Multichannel Speech Processing Handbook}, 2007.

\bibitem{serviere2004}
Christine Serviere and Dinh-Tuan Pham.
\newblock A novel method for permutation correction in frequency-domain in
  blind separation of speech mixtures.
\newblock In Carlos Puntonet and Alberto Prieto, editors, {\em Independent
  Component Analysis and Blind Signal Separation}, volume 3195 of {\em Lecture
  Notes in Computer Science}, pages 807--815. Springer Berlin / Heidelberg,
  2004.

\bibitem{1327173}
S.~Sanei, Wenwu Wang, and J.A. Chambers.
\newblock A coupled hmm for solving the permutation problem in frequency domain
  bss.
\newblock In {\em Acoustics, Speech, and Signal Processing, 2004. Proceedings.
  (ICASSP '04). IEEE International Conference on}, volume~5, pages V -- 565--8
  vol.5, may 2004.

\bibitem{wang:a}
Wenwu Wang, Jonathon~A. Chambers, and Saeid Sanei.
\newblock A novel hybrid approach to the permutation problem of frequency
  domain blind source separation.
\newblock In {\em ICA'04}, pages 532--539, 2004.

\bibitem{Donoho:1989aa}
David~L. Donoho and P.B. Stark.
\newblock Uncertainty principles and signal recovery.
\newblock {\em SIAM Journal on Applied Mathematics}, 49(3):906--931, 1989.

\bibitem{elad2002generalized}
M.~Elad and A.M. Bruckstein.
\newblock A generalized uncertainty principle and sparse representation in
  pairs of bases.
\newblock {\em Information Theory, IEEE Transactions on}, 48(9):2558--2567,
  2002.

\bibitem{Tao:2005aa}
Terence Tao.
\newblock An uncertainty principle for cyclic groups of prime order.
\newblock {\em Mathematical Research Letters}, 12:121--127, 2005.

\bibitem{Hall:1935aa}
Philip Hall.
\newblock On representatives of subsets.
\newblock {\em J. London Math. Soc.}, 10(1):26---30, 1935.

\bibitem{springerlink:10.1007/s00041-008-9042-0}
Joel Tropp.
\newblock On the linear independence of spikes and sines.
\newblock {\em Journal of Fourier Analysis and Applications}, 14:838--858,
  2008.
\newblock 10.1007/s00041-008-9042-0.

\bibitem{birkhoff1946three}
G.~Birkhoff.
\newblock Three observations on linear algebra, univ.
\newblock {\em Nac. Tucum{\'a}n. Rev. Ser. A}, 5:147--151, 1946.

\bibitem{von1953certain}
J.~Von~Neumann.
\newblock A certain zero-sum two-person game equivalent to the optimal
  assignment problem.
\newblock {\em Contributions to the Theory of Games}, 2:5--12, 1953.

\bibitem{berger1977geometrie}
M.~Berger.
\newblock {\em G{\'e}om{\'e}trie}, volume~4.
\newblock Cedic, 1977.

\bibitem{oxleymatroid}
J.~Oxley.
\newblock {\em Matroid theory}, volume~21.
\newblock Oxford University Press, USA, 1992.

\end{thebibliography}

\end{document}